\documentclass[12pt]{amsart}

\usepackage{amsmath,amssymb,amsthm}
\usepackage{url}
\usepackage{tikz}

\date{}
\author{Georg Braun}

\newtheorem{theorem}{Theorem}

\theoremstyle{definition}

\newtheorem{example}{Example}

\begin{document}

\title[Boolean percolation on digraphs]{Boolean percolation on digraphs and random exchange processes}

\thanks{{\em 2020 Mathematics Subject Classification.}  05C80, 60J80, 60K35, 82B43}  
\thanks{{\em Keywords: Boolean percolation; rumor spread and firework process; infinite paths in random graphs; long-range percolation; random exchange process; branching process with infinitely many types; recurrence and transience; spectral radius}}
\thanks{{\em Address:}
Mathematisches Institut,
Universit\"at T\"ubingen,
Auf der Morgenstelle 10,
72076 T\"ubingen, Germany.
Email: georg.braun@uni-tuebingen.de}

\begin{abstract}
We study, in a general graph-theoretic formulation, a long-range percolation model introduced by Lamperti in \cite{Lam70}. For various underlying directed graphs, we discuss connections between this model and random exchange processes. We clarify, for $n \in \mathbb{N}$, under which conditions the lattices $\mathbb{N}_0^n$ and $\mathbb{Z}^n$ are essentially covered in this model. Moreover, for all $n \geq 2$, we establish that it is impossible to cover the directed $n$-ary tree in our model.
\end{abstract}

\maketitle

\section{Introduction}

Percolation theory is a fascinating area of modern probability, which tries to understand under which conditions infinite components arise in random structures. In the present article, we study the properties of a Boolean percolation model on directed graphs and relate this model to a classical Markov chain known as the random exchange process.

Let $\mathcal{G} = (V,E)$ be a directed graph with an infinite, countable vertex set $V$. For all vertices $x$, $y \in V$, we denote by $d(x,y) \in \mathbb{N}_0 \cup \{ \infty \}$ the distance from $x$ to $y$ in $\mathcal{G}$. Note that $d: V \times V \rightarrow \mathbb{N}_0 \cup \{ \infty \}$ is an extended quasimetric on $V$, which is symmetric if and only if the graph $\mathcal{G}$ is undirected, i.e., $(x,y) \in E$ implies $(y,x) \in E$ for all $x$, $y \in V$. Moreover, for all $x \in V$ and $n \in \mathbb{N}_0$, we denote by $B_n (x)$ the open ball of radius $n$ starting from $x$, which is the set of all vertices $y \in V$ with $d(x,y) < n$.

Let $\mu = (\mu_n)_{n \in \mathbb{N}_0}$ be a probability vector and $(Y_x)_{x \in V}$ a family of i.i.d.\ random variables satisfying $\mathbb{P}[ Y_x = n] = \mu_n$ for all $n \geq 0$. In our percolation model, every vertex $x \in V$ will cover any vertex of $B_{Y_x} (x)$. Hence, the set of covered respectively uncovered vertices are
\[
V_\mu := V_\mu ( \mathcal{G})  := \bigcup\limits_{x \in V} B _{Y_x} (x)  \subseteq V, \quad \quad V_\mu^c := V_\mu^c ( \mathcal{G}) := V \setminus V_\mu ( \mathcal{G}).
\]
As we are interested in the properties of the random sets $V_\mu$ and $V_\mu^c$, we will always assume $\mu_0 \in (0,1)$, since $V_\mu = V$ almost surely in the case of $\mu_0=0$, and $V_\mu = \emptyset$ almost surely for $\mu_0=1$.

Let $x$, $y \in V$ and $V' = V_\mu$ or $V'= V_\mu^c$. Then, if both $x$ and $y$ are contained in $V'$ and connected by a path in $\mathcal{G}$, which uses only vertices from $V'$, we will say that $x$ and $y$ are in the same cluster.

To state our results, we introduce the following notation. Let $n \in \mathbb{N}$, $V = \mathbb{N}_0^n$ or $ V = \mathbb{Z}^n$, and $E$ be the set of all pairs $(x,x+e_j)$, where $x \in V$, $j=1,\ldots,n$, and $e_j = (\delta_{ij})_{i=1,\ldots,n}$. Then, we denote the resulting graph $\mathcal{G} = (V,E)$ by $\mathbb{N}_0^n$ respectively $\mathbb{Z}^n$. Furthermore, for all $n \geq 2$, we define the infinite directed $n$-ary tree $\mathcal{D}_n := (V_n, E_n)$ by
\begin{align*}
V_n &:= \bigcup\limits_{m \geq 0} \{ 1,\ldots,n \}^m, \quad \quad \quad \quad \quad \quad \quad \quad \quad \textnormal{where} \quad \hfill \{ 1,\ldots,n \}^0 := \emptyset,\\
E_n &:= \left\{ \left(\emptyset , 1 \right),\ldots, \left( \emptyset, n \right) \right\} \cup \left\{ \left( x, (x,j) \right)~\vert~ x \in V_n \setminus \{ \emptyset \}, ~j=1,\ldots,n \right\}.
\end{align*}

\begin{figure}[!h]
\centering
\begin{tikzpicture}
\draw (0,0) circle (4pt);
\draw (1,0) circle (4pt);
\draw (2,0) circle (4pt);
\draw (0,1) circle (4pt);
\draw (1,1) circle (4pt);
\draw (2,1) circle (4pt);
\draw (0,2) circle (4pt);
\draw (1,2) circle (4pt);
\draw (2,2) circle (4pt);
\draw[-stealth] (0.15,0)--(0.85,0);
\draw[-stealth] (1.15,0)--(1.85,0);
\draw[-stealth] (2.15,0)--(2.85,0);
\draw[-stealth] (0.15,1)--(0.85,1);
\draw[-stealth] (1.15,1)--(1.85,1);
\draw[-stealth] (2.15,1)--(2.85,1);
\draw[-stealth] (0.15,2)--(0.85,2);
\draw[-stealth] (1.15,2)--(1.85,2);
\draw[-stealth] (2.15,2)--(2.85,2);

\draw[-stealth] (0,0.15)--(0,0.85);
\draw[-stealth] (0,1.15)--(0,1.85);
\draw[-stealth] (0,2.15)--(0,2.85);
\draw[-stealth] (1,0.15)--(1,0.85);
\draw[-stealth] (1,1.15)--(1,1.85);
\draw[-stealth] (1,2.15)--(1,2.85);
\draw[-stealth] (2,0.15)--(2,0.85);
\draw[-stealth] (2,1.15)--(2,1.85);
\draw[-stealth] (2,2.15)--(2,2.85);

\begin{scope}[xshift=6cm]
\draw (0,0) circle (4pt);
\draw (1,0) circle (4pt);
\draw (2,0) circle (4pt);
\draw (0,1) circle (4pt);
\draw (1,1) circle (4pt);
\draw (2,1) circle (4pt);
\draw (0,2) circle (4pt);
\draw (1,2) circle (4pt);
\draw (2,2) circle (4pt);
\draw[-stealth] (-.65,0)--(-0.15,0);
\draw[-stealth] (0.15,0)--(0.85,0);
\draw[-stealth] (1.15,0)--(1.85,0);
\draw[-stealth] (2.15,0)--(2.85,0);
\draw[-stealth] (0.15,1)--(0.85,1);
\draw[-stealth] (-.65,1)--(-0.15,1);
\draw[-stealth] (1.15,1)--(1.85,1);
\draw[-stealth] (2.15,1)--(2.85,1);
\draw[-stealth] (0.15,2)--(0.85,2);
\draw[-stealth] (1.15,2)--(1.85,2);
\draw[-stealth] (2.15,2)--(2.85,2);
\draw[-stealth] (-.65,2)--(-0.15,2);
\draw[-stealth] (0,-0.65)--(0,-0.15);
\draw[-stealth] (1,-0.65)--(1,-0.15);
\draw[-stealth] (2,-0.65)--(2,-0.15);
\draw[-stealth] (0,0.15)--(0,0.85);
\draw[-stealth] (0,1.15)--(0,1.85);
\draw[-stealth] (0,2.15)--(0,2.85);
\draw[-stealth] (1,0.15)--(1,0.85);
\draw[-stealth] (1,1.15)--(1,1.85);
\draw[-stealth] (1,2.15)--(1,2.85);
\draw[-stealth] (2,0.15)--(2,0.85);
\draw[-stealth] (2,1.15)--(2,1.85);
\draw[-stealth] (2,2.15)--(2,2.85);
\end{scope}
\end{tikzpicture}
\caption[Illustration of $\mathbb{N}_0^2$ and $\mathbb{Z}^2$]{Illustration of the lattices $\mathbb{N}_0^2$ (left) and $\mathbb{Z}^2$ (right).}
\end{figure}
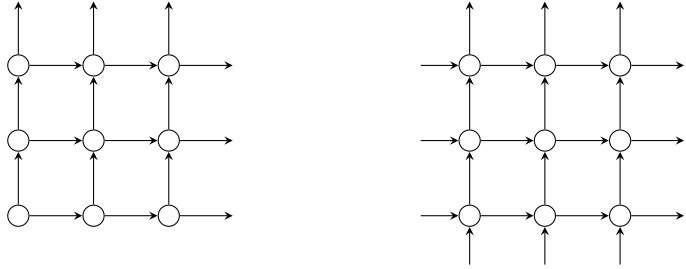

In this article, we will clarify under which conditions the graphs $\mathbb{N}_0^n$ and $\mathbb{Z}^n$ are (essentially) covered by a distribution $\mu$, compare Theorem 1 and Theorem 2 below. On the other hand, in Theorem 3, we will see that, for any distribution $\mu$ and $n \geq 2$, $\# V_\mu^c ( \mathcal{D}_n) = \infty$ almost surely.

To the best of our knowledge, the present percolation model was first studied by Lamperti in \cite{Lam70} for $\mathcal{G} = \mathbb{N}_0$. This research was motivated by statistical physics and included the following description. At each location $n \in \mathbb{N}_0$, there is a fountain, which sprays water to the right and is wetting the segment from $n+1$ to $n+ Y_n$. As $\mu_0> 0$, with some positive probability, a fountain fails to operate at all.

Our percolation model and variants of it were studied by various authors, compare \cite{LR08}, \cite{JMZ11}, \cite{JMZ14}, \cite{BZ13}, \cite{GGJR14}, and \cite{JMR19}. For a recent survey, see \cite{JMR19}. At this point, however, we want to postpone the discussion of how our new insights and results are related to these articles. 

We can interpret our percolation model as the spread of a rumor through a network, a firework process, or a discrete version of Boolean percolation. This model was introduced by Gilbert in \cite{Gil61}. First points are chosen randomly in $\mathbb{R}^n$ according to a Poisson point process. Then, in the simplest case, around these points, the unit sphere is covered. For monographs, which are concerned with Boolean percolation, see \cite{SKM87}, \cite{Hal88}, \cite{MR96}, and \cite{SW08}.

In the present article, we also investigate a connection between the above described graph-theoretic percolation model and a rather classical Markov chain, which is sometimes called random exchange process. As far as we know, it was first observed by Zerner in \cite[Section 1]{Zer18} that these two stochastic models are related to each other.

Let $(Y_n)_{n \geq 0}$ denote a sequence of i.i.d.\ random variables, which, as before, are distributed according to $\mu$. Then, we set $X_0 :=Y_0$ and recursively define
\[
X_{n+1} := \max \{ X_n - 1, Y_{n+1} \}, \quad n \in \mathbb{N}_0.
\]
To the best of our knowledge, this process $(X_n)_{n \geq 0}$ first occurred in a statistical research article on deepwater exchange of a fjord, see \cite{Gad73}. Later, it was studied, in more general form, in \cite{Hel74} and \cite{HN76}. In the following, we call the Markov chain $(X_n)_{n \geq 0}$ a (constant decrement) random exchange process. By construction, it has time-homogeneous transition probabilities and is irreducible on its state space $\mathcal{X}$, which is equal to $\mathbb{N}_0$ if $\mu$ is unbounded, and otherwise takes the form $\{ 0, 1, \ldots, n_0 \}$, where $n_0 := \sup \{n \in \mathbb{N}~\vert~\mu_n \neq 0 \}$. The transition matrix $P$ associated with $(X_n)_{n \geq 0}$ is
\[
P := P^\mu  := \left( P^\mu_{x,y} \right)_{x,y \in \mathcal{X}},~ \quad \textnormal{where} \quad P^\mu_{x,y} := \left\{ \begin{array}{l l}
\mu_y,&y \geq x,\\
\sum_{z=0}^{x-1} \mu_z, & y = x-1,\\
0,& y \leq x - 2. \\ 
\end{array} \right.
\]
As $(X_n)_{n \geq 0}$ respectively $P$ is irreducible, for all $z>0$, Green's function
\[
G(x,y \vert z) := \sum\limits_{n=0}^\infty P^n_{x,y}~z^n
\]
either converges or diverges simultaneously for all $x$, $y \in \mathcal{X}$, compare \cite[Chapter 1.1]{Woe00}. Therefore, independent of the choice of $x$, $y \in \mathcal{X}$, we can define the spectral radius of $(X_n)_{n \geq 0}$ respectively $P$ by
\[
\rho (P ) := \limsup\limits_{n \rightarrow \infty} \left( P^n_{x,y} \right)^{1/n} \in (0,1].
\]
More generally, if $A$ is an arbitrary irreducible matrix with nonnegative entries, we can define $\rho(A) \in [0,\infty]$ exactly in the same way.

Let us now state connections between the set $V_\mu$ of covered vertices in our percolation model and the Markov chain $(X_n)_{n \geq 0}$. We start by reformulating previous results in the following way.

\begin{theorem}
For any law $\mu$, the following statements are equivalent.\\ \\
\begin{tabular}{c l}
\textnormal{(i)}&Almost surely, $\# V_\mu^c ( \mathbb{Z}) < \infty$.\\ \\
\textnormal{(ii)}&Almost surely, $V_\mu ( \mathbb{Z}) = \mathbb{Z}$.\\ \\ 
\textnormal{(iii)}&The Markov chain $(X_n)_{n \geq 0}$ is not positive recurrent.\\ \\
\textnormal{(iv)}&The expectation of $\mu$ is infinite, i.e.,~~ $\sum_{n \geq 0} n \mu_n = \infty$.
\end{tabular}
\end{theorem}

It is not difficult to verify, more generally, that (i) and (ii) are equivalent if we replace $\mathbb{Z}$ by an arbitrary vertex-transitive graph.

By applying the Borel-Cantelli lemma, we can directly verify that (ii) and (iv) are equivalent statements. This equivalence was also observed, in a more general form, in \cite[Section 2.2]{JMZ11} and \cite[Section 4]{BZ13}. For any graph $\mathcal{G}=(V,E)$, we have $V_\mu = V$ almost surely if and only if
\[
\sum\limits_{x \in V} \sum\limits_{k \geq d(x,y)} \mu_k = \infty \quad \textnormal{for~all~} y \in V.
\]
For example, for all $n \geq 1$, $V_\mu ( \mathbb{Z}^n )= \mathbb{Z}^n$ almost surely if and only if the $n$-th moment of $\mu$ diverges. This kind of phenomenon is well-known in the context of Boolean percolation models, and, thus, it seems convenient to include some previous literature results at this point.

In \cite{Hal88}, Hall studied Boolean percolation on $\mathbb{R}^n$ with spheres of random i.i.d. radii and proved in \cite[Theorem 3.1]{Hal88} that the entire space is almost surely covered if and only if the $n$-th moment of the radius distribution diverges. In \cite{Gou08}, Gour established that if the $n$-th moment of the radius distribution is finite, there exists a critical value for the intensity of the underlying Poisson process. Recently, more results on phase transitions were deduced in \cite{ATT18} and \cite{DRT20}. However, there are also results on other aspects of Boolean percolation. For example, in \cite{ABGM14}, it was shown that this model is noise sensitive, and in \cite{LPZ17}, the capacity functional was studied.

In \cite{ARS04} and \cite{Bez21}, Boolean percolation was studied on $[0,\infty)^n$ when, instead of the sphere around a point $x$, the set $x + [0,R_x)^d$ is occupied, where $R_x$ is the radius associated which $x$. The results in \cite{ARS04} characterize under which conditions the entire space is essentially covered, and interestingly depend on whether $n=1$ or $n \geq 2$. In \cite{Bez21}, Bezborodov observed, for $n=1$ and some radius distributions, that the covered volume fraction is one, but all clusters are bounded almost surely.

In \cite{CG14}, Coletti and Grynberg studied a model on $\mathbb{Z}^n$, in which first Bernoulli percolation with parameter $p \in (0,1)$ is performed, and then, independently, around the present points, random i.i.d.\ balls are covered. Again, the occupied region is almost surely $\mathbb{Z}^n$ if and only if the $n$-th moment of the radius distribution diverges. For a study of this percolation model on doubling graphs, also see \cite{CGM20}.

Let us return to Theorem 1. The equivalence of (iii) and (iv) was first observed by Helland in \cite[Section 3]{Hel74} and also mentioned by Kellerer in \cite[comments after Theorem 2.6]{Kel06}. We can deduce it as follows. Due to the form of the transition probabilities of the Markov chain $(X_n)_{n \geq 0}$, any invariant measure $\tau = (\tau_x)_{x \in \mathcal{X}}$ has to satisfy
\[
\tau_{x} = \sum\limits_{z=0}^x \tau_z \mu_x  + \tau_{x+1} \sum\limits_{z=0}^{x} \mu_z, \quad \textnormal{provided~that~$x$,~$x+1$} \in \mathcal{X}.
\]
Solving this recurrence relation yields the representation
\begin{align}
\tau_x = \tau_0 \left( \sum\limits_{z \geq x} \mu_z \right) \left( \prod\limits_{y=0}^{x-1} \sum\limits_{z=0}^y \mu_z \right)^{-1}, \quad x \in \mathcal{X}.
\end{align}
By a careful look at this formula indeed, it follows that (iii) and (iv) are equivalent. Moreover, if the distribution $\mu$ has a finite expectation, we can determine the stationary solution of $(X_n)_{n \geq 0}$ from (1) via normalization. In Section 3, we will present some concrete examples.

For results on positive recurrence of more general exchange processes with random decrements, the reader may consult \cite[Section 2]{HN76}.

\begin{theorem}
For any law $\mu$, the following statements are equivalent.\\ \\
\begin{tabular}{c l}
\textnormal{(a)}&There exists $n \in \mathbb{N}$ with $\# V_\mu^c ( \mathbb{N}_0^n) < \infty$ almost surely.\\ \\
\textnormal{(b)}&For all $n \in \mathbb{N}$, $\# V_\mu^c ( \mathbb{N}_0^n) < \infty$ almost surely.\\ \\
\textnormal{(c)}&The Markov chain $(X_n)_{n \geq 0}$ is transient.\\ \\
\textnormal{(d)}&$\sum_{m \geq 0} \prod_{k=1}^m \sum_{l=0}^{k-1} \mu_l < \infty$.\\ \\
\end{tabular}\\
Moreover, if one of these conditions is satisfied, then $\mathbb{E} [\# V_\mu^c ( \mathbb{N}_0^n)] < \infty$ for all $n \in \mathbb{N}$ and there exists $\alpha \in (0,\infty)$ with $\mathbb{E} [ \exp ( \alpha \# V_\mu^c ( \mathbb{N}_0))] < \infty$.
\end{theorem}

This theorem improves on previous works by revealing that, rather surprisingly, the value of $n \in \mathbb{N}$ does not influence whether all but finitely many points of the graph $\mathbb{N}_0^n$ are covered by a distribution $\mu$. In the appendix of \cite{Lam70}, Kesten proved that $\# V_\mu^c ( \mathbb{N}_0) < \infty$ almost surely if and only if condition (d) in Theorem 2 is satisfied. This result was later rediscovered by various authors, partly in a different and more general form, compare \cite[comments to Proposition 6.6]{Kel06}, \cite[Theorem 2.1]{JMZ11}, \cite[Theorem 1]{GGJR14}, and \cite[Section 3]{BZ13}.

As observed by Zerner in \cite[Proposition 1.1]{Zer18} and suggested by our notation, we can couple the set of covered points $V_\mu ( \mathbb{N}_0)$ and the random exchange process $(X_n)_{n \geq 0}$ by using the same sequence of random variables $(Y_n)_{n \geq 0}$ in both definitions. Then, by construction,
\begin{center}
\begin{tabular}{l l l}
$V_\mu ( \mathbb{N}_0)$&$=$&$ \left\{ n \in \mathbb{N}_0~\vert~ \exists k \in \{0,\ldots,n\} : Y_k > n - k \right\}$\\[5pt]
&$=$&$ \bigl\{ n \in \mathbb{N}_0~\bigm\vert~ \max_{0 \leq k \leq n}  \left( Y_k - (n-k) \right) > 0 \bigr\}$\\[5pt]
&$=$&$ \{ n \in \mathbb{N}_0~\vert~ X_{n} > 0 \}$.
\end{tabular}
\end{center}
Consequently, we know that the Markov chain $(X_n)_{n \geq 0}$ is transient if and only if $\# V_\mu^c ( \mathbb{N}_0 ) < \infty$ almost surely.

Let $n \geq 2$ and consider the infinite directed $n$-ary tree $\mathcal{D}_n = (V_n, E_n)$. Then, interestingly, we can associate a multitype branching process $(Z_m)_{m \geq 0}$ to our percolation model on $\mathcal{D}_n$ in the following way.

Let $k \in \mathbb{N}_0$ and $y \in V_n$ with $d(\emptyset,y)=k$. Then, we identify the vertex $y$ with an individual of the $k$-th generation of $(Z_m)_{m \geq 0}$ if and only if $y \in V_\mu$ and $y$ is contained in the same cluster as $\emptyset$. In other words, we demand that all vertices, which form the path from $\emptyset$ to $y$ in $\mathcal{D}_n$, are contained in $V_\mu ( \mathcal{D}_n)$. If this condition is satisfied, we define the type of $y$ by
\[
z_y:= \max \left\{ Y_x - d(x,y)  ~\vert~ x \in V_n,~d(x,y) < \infty \right\}.
\]
By construction, if $Y_\emptyset \geq 1$, the branching process $(Z_n)_{n \geq 0}$ starts with one individual of type $Y_\emptyset$. However, on the event $Y_\emptyset = 0$, there are no individuals at all. For a vertex $y \in V_n$ to be identified with an individual in $(Z_m)_{n \geq 0}$, necessarily $y \in V_\mu$, i.e., there exists $x \in V_n$ with $Y_x > d(x,y)$. Hence, as $\mu_0 \in (0,1)$, the type space of the branching process $(Z_m)_{n \geq 0}$ is $\mathcal{Z} = \mathcal{X} \setminus \{ 0 \}$.

We can describe the reproduction in this branching process as follows. Every individual of type $x \geq 2$ has exactly $n$ children, whose types are independent of each other. For each of them, the probability of type $y \in \mathcal{Z}$ is $M_{x,y}:= P_{x,y}$. On the other hand, an individual of type $1$ has $n$ potential children, which are again independent of each other. For all $z \in \mathcal{Z}$, the probability that a given potential child is born and of type $z$ is $M_{1,z}:=P_{1,z} $. However, with probability $\mu_0$, a potential child is not born.

As the type space $\mathcal{Z}$ is infinite in general, we distinguish between the local and global extinction of $(Z_m)_{m \geq 0}$. This process dies out globally if, at some moment, the total number of individuals vanishes. It dies out locally if, for all $z \in \mathcal{Z}$, only finitely many individuals of type $z$ are born. While global extinction always implies local extinction, the reverse is not true in general for branching processes with infinitely many types.

\begin{theorem}
Let $n \geq 2$. Then, for any distribution $\mu$,  $\# V_\mu^c ( \mathcal{D}_n) = \infty$ almost surely. Moreover, the following statements are equivalent.\\ \\
\begin{tabular}{c l}
\textnormal{(A)}&Almost surely, $V_\mu ( \mathcal{D}_n)$ contains a path of infinite length.\\ \\
\textnormal{(B)}&With a positive probability, $(Z_m)_{m \geq 0}$ will not die out globally.\\ \\
\textnormal{(C)}&With a positive probability, $(Z_m)_{m \geq 0}$ will not die out locally.\\ \\
\textnormal{(D)}&$\rho (M) > n^{-1}$, where $M:= (M_{x,y})_{x,y \in \mathcal{Z}}$.
\end{tabular}
\end{theorem}

Note that, up to multiplication with $n \geq 2$, $M$ is the mean matrix of the branching process $(Z_m)_{m \geq 0}$. To some degree, this explains why condition (D) is related to (B) and (C). Also, observe that $M$ arises from the transition matrix $P$ of the random exchange process $(X_n)_{n \geq 0}$ simply by deleting both the first row and column.

For an introduction to infinite type branching processes, we recommend Braunsteins' exposition in \cite[Chapter 2]{Bra18} and the references mentioned therein. This presentation also explains the rather well-understood results on the extinction of finite type branching processes.

We also want to note that the branching processes $(Z_n)_{n \geq 0}$, which we consider in the present article, have a mean matrix of upper Hessenberg form. Recently, Braunsteins and Haupthenne studied the extinction of branching processes with a lower Hessenberg mean matrix in \cite{BH19}.

\section{Proof of Theorem 2 and Theorem 3}

\begin{proof}[Proof of Theorem 2]
(b)$\Longrightarrow$(a). This implication is clear.

(b)$\Longrightarrow$(c)$\Longrightarrow$(a). From our coupling between $V_\mu ( \mathbb{N}_0)$ and $(X_n)_{n \geq 0}$, we know that $(X_n)_{n \geq 0}$ is transient if and only if $V_\mu^c ( \mathbb{N}_0) < \infty$ almost surely. In particular, the implications (b)$\Longrightarrow$(c) and (c)$\Longrightarrow$(a) follow.

(d)$\Longrightarrow$(a). By Kolmogorov's 0-1 law, for any probability distribution $\mu$, either $\# V_\mu^c ( \mathbb{N}_0) < \infty$ almost surely or $\#V_\mu^c ( \mathbb{N}_0) = \infty$ almost surely. By identifying the expression in (d) with $\mathbb{E} [\# V_\mu^c( \mathbb{N}_0)]$, the implication follows.

Finally, let us assume that condition (a) holds for a distribution $\mu$.

In the first step, we verify that we can restrict ourselves to the case $n=1$. For this, suppose $\# V_\mu^c ( \mathbb{N}_0^n) < \infty$ almost surely for some $n \geq 2$. Then, consider the subgraph $\mathcal{G}'=(V',E')$ of $\mathbb{N}_0^n$, which is induced by the vertex set $V'$ of all $(x_1,\ldots,x_n) \in \mathbb{N}_0^n$ with $x_j = 0$ for all $j=2,\ldots,n$. By construction, $\mathcal{G}'$ is isomorphic to $\mathbb{N}_0$, and we know that $V_\mu^c ( \mathcal{G}') = V_\mu^c ( \mathbb{N}_0^n) \cap V'$ is finite almost surely. Consequently, $\# V_\mu^c ( \mathbb{N}_0 ) < \infty$ almost surely.

In the second step, we prove all remaining claims. As $0 < \mu_0 < 1$,
\[
p:= \mathbb{P} \left[ V_\mu(\mathbb{N}_0)=\mathbb{N} \right] > 0.
\]
Therefore, by the strong Markov property, we know that $\# V_\mu^c ( \mathbb{N}_0)$ is geometrically distributed with parameter $p$. In particular, this random variable has a finite exponential moment, and we also deduce (d), i.e.,
\begin{align}
\mathbb{E} \left[ \# V_\mu^c (\mathbb{N}_0) \right] = \sum\limits_{m \geq 0} q_m < \infty, \quad \textnormal{where} \quad q_m := \mathbb{P }\left[m \in V_\mu^c ( \mathbb{N}_0) \right]. 
\end{align}
Let $n \in \mathbb{N}$ and $x=(x_1,\ldots,x_n) \in \mathbb{N}_0^n$. For all $j=1,\ldots,n$, let $\pi_j$ denote the unique path from $(x_1,\ldots,x_{j-1},0,x_{j+1},\ldots,x_n)$ to $x$. Then, for all $j=1,\ldots,n$, the path $\pi_j$ consists of $x_j$ edges in direction $e_j$, and, for all $i \neq j$, the paths $\pi_i$ and $\pi_j$ only share one vertex, which is their endpoint $x =(x_1,\ldots,x_n)$. So, for the event $\{ x \in V_\mu^c (\mathbb{N}_0^n) \}$ to occur, it is necessary that for each $j=1,\ldots,n$ there exists no vertex $y$ contained in the path $\pi_j$ such that $Z_y > d(y,x)$. This defines $n$ independent events, whose probabilities can be described with the sequence $(q_m)_{m \geq 0}$ defined in (2). All in all,\\
\begin{tabular}{l l l}
$\displaystyle \mathbb{E} \left[ \# V_\mu^c ( \mathbb{N}_0^n) \right]$&$=$&$ \displaystyle \sum\limits_{x \in \mathbb{N}_0^n} \mathbb{P} \left[ x \in V_\mu^c ( \mathbb{N}_0^n) \right] \leq \sum\limits_{(x_1,\ldots,x_n) \in \mathbb{N}_0^n}~ \prod\limits_{j=1}^n~ q_{x_j}$\\ \\
&$=$&$\displaystyle \sum\limits_{x_1 \geq 0} q_{x_1} \sum\limits_{x_2 \geq 0} q_{x_2} \cdots \sum\limits_{x_{n-1} \geq 0} q_{x_{n-1}} \sum\limits_{x_n \geq 0} q_{x_n}$\\ \\
&$=$&$\displaystyle \mathbb{E} \left[\# V_\mu^c ( \mathbb{N}_0) \right]^n < \infty $.\\~\\
\end{tabular}~\\
In particular, for all $n \geq 2$, $\# V_\mu^c ( \mathbb{N}_0^n) < \infty$ almost surely, and condition (b) holds. Since we have already verified (b)$\Longrightarrow$(c), this finishes the proof.
\end{proof}

\begin{proof}[Proof of Theorem 3]
First, we verify that for all $n \geq 2$ and any law $\mu$, $\# V_\mu^c ( \mathcal{D}_n) = \infty$ almost surely. For this, for all $m \in \mathbb{N}$, we set
\[
r_m := \mathbb{P} \left[ \exists y \in V_\mu^c ( \mathcal{D}_n):~d(\emptyset,y)=m  \right].
\]
As $0 < \mu_0 < 1$, we know $r_m \in (0,1)$ for all $m \in \mathbb{N}$. Moreover, for all $j=1,\ldots,N$, we denote by $\mathcal{G}_j$ the induced subgraph obtained from $\mathcal{D}_n$ by restricting to all vertices, which can be reached from $j \in V_n$.

Let $m \geq 2$. Then, we know that there exists a $y \in V_\mu^c ( \mathcal{D}_n)$ with $d(\emptyset,y)=m$ if and only if $Y_\emptyset \leq m$ and, for some $j=1,\ldots,n$, there exists a vertex $z$ in the graph $\mathcal{G}_j$ with $d(j,z)=m-1$, which is not covered by any other vertex of $\mathcal{G}_j$. Note that the latter event is independent of $Y_\emptyset$ and that the graphs $\mathcal{G}_1,\ldots,\mathcal{G}_n$ are isomorphic to $\mathcal{D}_n$. Consequently, for all $m \geq 1$, we obtain the recurrence relation
\begin{equation}
r_{m+1} = \left( 1 - (1 - r_m)^n \right)~ F(m+1), \quad \quad \textnormal{where} \quad F(k) := \sum\limits_{l=0}^k \mu_l.
\end{equation}
Let $N \in \mathbb{N}$ with $F(N) > 1/2$. Then, by (3), for all $m \geq N$,
\begin{align}
r_{m+1} \geq \left( 1 -(1-r_m)^2 \right) F(N) = r_m (2-r_m) F(N). 
\end{align}
The map $f_N: [0,1] \rightarrow [0,1]$, $x \mapsto x (2-x) F(N)$, is monotone increasing. Hence, due to the estimate (4), iteration of the function $f_N$ yields
\[
r_m \geq f_N^{m-N} (r_N) \quad \textnormal{for~all~} m \geq N+1.
\]
The map $f_N$ has the two fixpoints $0$ and $x_N:= 2- F(N)^{-1} \in (0,1]$. Hence, by monotonicity, if $r_N \geq x_N$, then also $r_m \geq x_N$ for all $m \geq N$. On the other hand, if $r_N < x_N$, then, since $f_N$ is concave and $f'_N (0) > 1$, $f_N^k ( r_N) \rightarrow x_N$ for $k \rightarrow \infty$. In both cases, we can deduce
\[
\liminf\limits_{m \rightarrow \infty} r_m \geq x_N = 2- F(N)^{-1}.
\]
As $N \in \mathbb{N}$ can be chosen arbitrarily large in this argument, it follows that $r_m \rightarrow 1$ for $m \rightarrow \infty$. In particular, $V_\mu^c ( \mathcal{D}_n)$ is almost surely non-empty.

By Kolmogorov's 0-1 law, we know that either $\# V_\mu^c ( \mathcal{D}_n)$ is finite almost surely, or this random variable is infinite almost surely. In the first case, due to our assumption $\mu_0 \in (0,1)$, it would follow that $V_\mu^c ( \mathcal{D}_n)$ is empty with a positive probability. Hence, we can conclude $\# V_\mu^c ( \mathcal{D}_n) = \infty$ almost surely.\\

In the second step of this proof, we now verify that indeed the statements (A), (B), (C), and (D) are equivalent to each other.

(C)$\Longrightarrow$(B). This implication is clear.

(A)$\Longleftrightarrow$(B). If (A) holds, then, with a positive probability, $V_\mu( \mathcal{D}_n)$ contains an infinite path starting from the root $\emptyset$. On this event, $(Z_n)_{n \geq 0}$ does not die out globally, i.e., condition (B) holds. Conversely, if (B) holds, then, with a positive probability, $V_\mu ( \mathcal{D}_n)$ contains an infinite path. By Kolmogorov's $0$-$1$ law, (A) follows.

(C)$\Longleftrightarrow$(D). Since the mean matrix $M$ of the branching process $(Z_n)_{n \geq 0}$ is irreducible, this equivalence follows from the theory of multitype branching processes, compare \cite[Theorem 9]{Bra18}, \cite{GM06}, or \cite{BZ09}.

(B)$\Longrightarrow$(C). Suppose, for some distribution $\mu$, that $(Z_n)_{n \geq 0}$ dies out locally almost surely but survives forever with a positive probability. Then, as $(Z_n)_{n \geq 0}$ starts with a single individual of random type $Y_\emptyset$, with some positive probability, $(Z_n)_{n \geq 0}$ survives forever, and no individuals of type $1$ are born. On this event, we would know that $\# V_\mu^c ( \mathcal{D}_n) < \infty$, and this is a contradiction to our first claim. The implication follows.
\end{proof}

\section{Examples}

\begin{example}
Let $m \in \mathbb{N}$, $m \geq 2$, and $\mu$ be the uniform distribution on $\{ 0,1,\ldots,m-1\}$. Then, by (1), the stationary solution $\tau$ of $(X_n)_{n \geq 0}$ is
\[
\tau_n = \frac{m!}{m^m} (m-n)  \frac{m^{n-1}}{n!}, \quad n \in \{ 0, \ldots, m-1 \}.
\]
This law $\tau$ is a terminating member of the Kemp family of generalized hypergeometric probability distributions, compare \cite[Section 2.4.1]{JKK05}. However, it also naturally arises from Naor's urn model \cite[Appendix]{Nao57}, also see \cite[Section 11.2.12]{JKK05}. Assume that there are $m$ balls in an urn, of which one is red, and the rest are white. In each step, pick one ball, and if it is white, replace it with a red ball. Continue until the first time $T$, at which a red ball gets chosen. Then, the distribution of $T$ is
\[
\mathbb{P} [T=n] = (m-1)!~m^{-n}~\frac{n}{(m-n)!}, \quad n \in \{ 1,\ldots,m \},
\]
and $m-T$, i.e.\ the number of tries not needed, has distribution $\tau$.
\end{example}

\begin{example}
Let $p \in (0,1)$ and $\mu$ be the geometric distribution with parameter $1-p$. Then, by (1), the stationary solution $\tau$ of $(X_n)_{n \geq 0}$ is
\[
\tau_n = \tau_0~ p^n \Bigg( \prod_{k=1}^{n} \left( 1 - p^{k} \right)  \Bigg)^{-1} = \tau_0 ~ \frac{p^n}{(p;p)_n}, \quad n \in \mathbb{N}_0,
\]
where $(a;q)_n$ is the $q$-Pochhammer symbol. By normalisation,
\[
\tau_0 = \Bigg( \sum\limits_{n \geq 0} \frac{p^n}{(p;p)_n} \Bigg)^{-1} = \frac{1}{(p;p)_{\infty}} = \phi(p)^{-1},
\]
where we have applied the $q$-binomial theorem, and $\phi$ denotes Euler's function. In \cite[Section 4]{BB88}, Benkherouf and Bather discussed, in more general form, this distribution $\tau$ and referred to it as an Euler distribution. For more information, also see \cite[Section 10.8.2]{JKK05}.
\end{example}

\begin{example}
Assume that there exists $c \in (0,\infty)$ and $n_0 \in \mathbb{N}$ with
\[
\sum\limits_{k>n} \mu_k = \frac{c}{n} \quad \textnormal{for~all~} n \geq n_0.
\]
Then, as the expectation of $\mu$ is infinite, we know that the statements (i)-(iv) from Theorem 1 do hold. Moreover, it follows from the Gaussian ratio test that the condition (d) in Theorem 2 is satisfied if and only if $c>1$. According to \cite[Theorem 3.2]{HN76}, for any value of $c \in (0,\infty)$,
\[
\lim\limits_{n \rightarrow \infty} \mathbb{P} [ X_n~n^{-1} \leq y] = y^c (y+1)^{-c} \quad \textnormal{for~all~} y \in (0,\infty).
\]
This limit is an inverse Beta distribution with $\alpha=c$ and $\beta=1$.
\end{example}

\begin{example}
Assume that $\mu = (\mu_n)_{n \geq 0}$ has finite support. Then, $(X_n)_{n \geq 0}$ has a stationary solution and $\rho(P)=1$. Moreover, $\rho(M)$ is the spectral radius and Perron Frobenius eigenvalue of $M$. As in the proof of Theorem 3, we define, for the case $n=2$,
\[
r_m:= \mathbb{P} \left[ \exists y \in V_\mu^c ( \mathcal{D}_2): d(\emptyset,y)=m \right], \quad m \in \mathbb{N}.
\]
Then, observe that, for all $m \geq n_0 := \sup \{ n \in \mathbb{N}~\vert~ \mu_n \neq 0 \}$, the recurrence relation (3) simplifies into
\[
r_{m+1}  = ( 1- (1-r_m)^2 ) =  r_m ( 2 -r_m).
\]
This recursion is a modified version of the logistic equation. It follows that
\[
r_m = 1 - \exp \left( - c~ 2^m \right) \quad \textnormal{for~all~} m \geq n_0,
\]
where $c \in (0,\infty)$ is a fixed parameter.
\end{example}

\begin{example}
Let $n \in \mathbb{N}$, $p \in (0,1)$, $\mu_n := p$ and $\mu_0:=1-p$. Then, the matrix $M=M_{n,p}$ has dimension $n$ and is of the form
\[
M_{n,p} = \begin{pmatrix}
0&0&\cdots&\cdots&0&p\\
1-p&0&\cdots&\cdots&0&p\\
0&1-p&0&\cdots&0&p\\
\vdots&&\ddots&&&\vdots\\
\vdots&&&\ddots&&\vdots\\
0&\cdots&\cdots&0&1-p&p
\end{pmatrix}
\]
The characteristic polynomial $\chi_{n,p} =\chi_{n,p} (z)$ of $M_{n,p}$ satisfies
\[
\chi_{n,p} (z) = \det ( z \mathbf{1}_{n} - M_{n,p}) = z \chi_{n-1,p} (z) -p (1-p)^{n-1}, 
\]
and this recurrence relation can be deduced from a Laplace expansion of the first row of $M_{n,p}$. It follows from $\xi_{1,p}(z)=z-p$, that
\[
\chi_{n,p} (z) = \frac{p (1-p)^n + (z-1) z^n}{p+z-1}.
\]
We know that $\rho ( M_{n,p})$ is the largest zero of this polynomial in $(0,1)$.
\end{example}

\noindent \textbf{Acknowledgment.} The author thanks Martin Zerner and Elmar Teufl for many helpful comments. The author is grateful for a PhD scholarship of the Landesgraduiertenf\"orderung Baden-W\"urttemberg.\\



\begin{thebibliography}{99}
\bibitem{ABGM14} \textsc{Daniel Ahlberg, Erik I. Broman, Simon Griffiths, and Robert Morris.} (2014). Noise sensitivity in continuum percolation. \textit{Israel J. Math.}, 201(2), 847--899. MR3265306.

\bibitem{ATT18} \textsc{Daniel Ahlberg, Vincent Tassion, and Augusto Q. Teixeira} (2018). Sharpness of the phase transition for continuum percolation in $\mathbb {R}^ 2$. \textit{Probab. Theory Related Fields}, 172(1-2), 525--581. MR3851838.

\bibitem{ARS04} \textsc{Siva R. Athreya, Rahul Roy, and Anish Sarkar} (2004). On the Coverage of Space by Random Sets. \textit{Adv. in Appl. Probab.} 36(1), 1--18. MR2035771.

\bibitem{BB88} \textsc{Lakdere Benkherouf, and John A. Bather} (1988). Oil Exploration: Sequential Decisions in the Face of Uncertainty. \textit{J. Appl. Probab.}, 25(3), 529--543. MR954501.

\bibitem{BZ09} \textsc{Daniela Bertacchi, and Fabio Zucca} (2009). Characterization of Critical Values of Branching Random Walks on Weighted Graphs through Infinite-Type Branching Processes. \textit{J. Stat. Phys.}, 134(1), 53--65. MR2489494.

\bibitem{BZ13} \textsc{Daniela Bertacchi, and Fabio Zucca.} (2013). Rumor Processes in Random Environment on $\mathbb{N}$ and on Galton-Watson Trees. \textit{J. Stat. Phys.}, 153(3), 486--511. MR3107655.

\bibitem{Bez21} \textsc{Viktor Bezborodov} (2021). Non-triviality in a totally asymmetric one-dimensional Boolean percolation model on a half-line. \textit{Statist. Probab. Lett.},  176, Paper No. 109155. MR4263133.

\bibitem{Bra18} \textsc{Peter T. Braunsteins} (2018). \textit{Extinction in branching processes with countably many types.} PhD thesis, University of Melbourne. Available at \url{https://minerva-access.unimelb.edu.au/handle/11343/210538}

\bibitem{BH19} \textsc{Peter T. Braunsteins, and Sophie Haupthenne} (2019). Extinction in lower Hessenberg processes with countably many types. \textit{Ann. Appl. Probab.}, 29(5), 2782--2818. MR4019875.

\bibitem{CG14} \textsc{Christian F. Coletti, and Sebastian P. Grynberg} (2014). Absence of percolation in the Bernoulli Boolean model. \textit{Preprint}. Available at \url{https://arxiv.org/abs/1402.3118}.

\bibitem{CGM20} \textsc{Christian F. Coletti, Daniel Miranda, and Sebastian P. Grynberg} (2020). Boolean Percolation on Doubling Graphs. \textit{J. Stat. Phys.}, 178(3), 814--831. MR4059962.

\bibitem{DRT20} \textsc{Hugo Duminil-Copin, Aran Raoufi, and Vincent Tassion.} (2020). Subcritical phase of $d$-dimensional Poisson-Boolean percolation and its vacant set. \textit{Ann. H. Lebesgue} 3, 677--700. MR4149823.

\bibitem{Gad73} \textsc{Herman G. Gade} (1973). Deep Water Exchanges in a Sill Fjord: A Stochastic Process. \textit{J. Phys. Occanogr.}, 3, 213--219.

\bibitem{GGJR14} \textsc{Sandro Gallo, Nancy L. Garcia, Valdivino V. Junior, and Pablo M. Rodr\'{\i}gues} (2014). Rumor processes on $\mathbb{N}$ and Discrete Renewal Processes. \textit{J. Stat. Phys.}, 155(3), 591--602. MR3192175.

\bibitem{GM06} \textsc{Nina Gantert, and Sebastian M{\"u}ller} (2006). The critical Markov chain is transient. \textit{Markov Process. Related Fields}, 12(4), 805--814. MR2284404.

\bibitem{Gil61} \textsc{Edgar N. Gilbert} (1961). Random Plane Networks. \textit{J. Soc. Indust. Appl. Math.}, 9, 533--543. MR132566.

\bibitem{Gou08} \textsc{Jean-Baptiste Gou\'{e}r\'{e}} (2008). Subcritical regimes in the Poisson Boolean model of continuum percolation. \textit{Ann. Probab.}, 36(4), 1209--1220. MR2435847.

\bibitem{Hal85} \textsc{Peter Hall} (1985). On Continuum Percolation. \textit{Ann. Probab.}, 13(4), 1250--1266. MR806222.

\bibitem{Hal88} \textsc{Peter Hall} (1988). \textit{Introduction to the Theory of Coverage Processes}. Wiley Series in Probability and Mathematical Statistics: Probability and Mathematical Statistics. John Wiley \& Sons, Inc., New York. MR973404.

\bibitem{Hel74} \textsc{Inge S. Helland} (1974). A random exchange model with constant decrements. \textit{Report}, University of Bergen. Available at \url{https://bora.uib.no/bora-xmlui/handle/1956/19470}.

\bibitem{HN76} \textsc{Inge S. Helland, and Trygve S. Nilsen} (1976). On a general random exchange model. \textit{J. Appl. Probability}, 13(4), 781--790. MR431437.

\bibitem{JKK05} \textsc{Norman L. Johnson, Adrienne W. Kemp, and Samuel Kotz} (2005). \textit{Univariate Discrete Distributions}. Wiley Series in Probability and Statistics. Third edition. Wiley-Interscience, John Wiley \& Sons, Hoboken, NJ. MR2163227.

\bibitem{JMR19} \textsc{Valdivino V. Junior, F\'{a}bio P. Machado, and Krishnamurthi Ravishankar} (2019). The Rumor Percolation Model and Its Variations. In: \textit{Sojourns in Probability Theory and Statistical Physics - II} (edited by V. Sidoravicius), pages 208-227. Springer Proceedings in Mathematics \& Statistic. Springer, Singapore.

\bibitem{JMZ11} \textsc{Valdivino V. Junior, F\'{a}bio P. Machado, and Mauricio Zuluaga} (2011). Rumor processes on $\mathbb{N}$. \textit{J. Appl. Probab.}, 48(3), 624--636. MR2884804.

\bibitem{JMZ14}  \textsc{Valdivino V. Junior, F\'{a}bio P. Machado, and Mauricio Zuluaga} (2014). The cone percolation on $\mathbb{T}_d$. \textit{Braz. J. Probab. Stat.}, 28(3), 367--375. MR3263053.

\bibitem{Kel06} \textsc{Hans G. Kellerer} (2006). Random Dynamical Systems on Ordered Topological Spaces. \textit{Stoch. Dyn.}, 6(3), 255--300. MR2258486.

\bibitem{Lam70} \textsc{John W. Lamperti} (1970). Maximal Branching Processes and `Long-Range Percolation'. \textit{J. Appl. Probability}, 7(1), 89--98. MR254930.

\bibitem{LPZ17} \textsc{G\"{u}nter Last, Mathew D. Penrose, and Sergei Zuyev} (2017). On the capacity functional of the infinite cluster of a Boolean model. \textit{Ann. Appl. Probab.}, 27(3), 1678--1701. MR3678482.

\bibitem{LR08} \textsc{\'{E}lcio Lebensztayn, and Pablo M. Rodr\'{\i}gues} (2008). The disk-percolation model on graphs. \textit{Statist. Probab. Lett.}, 78(14), 2130--2136. MR2458022.

\bibitem{MR96} \textsc{Ronald Meester, and Rahul Roy} (1996). \textit{Continuum Percolation.} Volume 119 of Cambridge Tracts in Mathematics. Cambridge University Press, Cambridge.

\bibitem{Nao57} \textsc{Pinhas Naor} (1957). Normal Approximation to Machine Interference with Many Repairmem. \textit{J. R. Statist. Soc. Ser. B}, 19(2), 334--341.

\bibitem{PY01} \textsc{Mathew D. Penrose, and Joseph E. Yukich} (2002). Limit Theory for Random Sequential Packing and Deposition. \textit{Ann. Appl. Probab.}, 12(1), 272--301. MR1890065.

\bibitem{SW08} \textsc{Rolf Schneider, and Wolfgang Weil} (2008). \textit{Stochastic and Integral Geometry.} Probability and Its Applications. Springer-Verlag, Berlin. MR2455326.

\bibitem{Sep00} \textsc{Timo Sepp\"{a}l\"{a}inen} (2000). Strong law of large numbers for the interface in ballistic deposition. \textit{Ann. Inst. H. Poincar\'{e} Probab. Statist.}, 36(6), 691--736. MR1797390.

\bibitem{SKM87} \textsc{Dietrich Stoyan, Wilfrid S. Kendall, and Joseph Mecke} (1987). \textit{Stochastic Geometry and its Applications}. Wiley Series in Probability and Mathematical Statistics: Applied Probability and Statistics. John Wiley \& Sons, Ltd., Chichester. MR895588.

\bibitem{Woe00} \textsc{Wolfgang Woess} (2000). \textit{Random Walks on Infinite Graphs and Groups.} Cambridge Tracts in Mathematics. Cambridge University Press, Cambridge. MR1743100.

\bibitem{Zer18} \textsc{Martin P. W. Zerner} (2018). Recurrence and transience of contractive autoregressive processes and related Markov chains. \textit{Electron. J. Probab.}, 23(27), 1--24. MR3779820.
\end{thebibliography}
\end{document}